\documentclass{amsart}

\input{xypic.tex} \xyoption{all}
\usepackage{amsmath,amssymb,amsfonts,amstext,amsthm,latexsym}

\newcommand{\nwc}{\newcommand}

\nwc{\cin}{\textbf{(v)}}

\nwc{\doi}{\textbf{(ii)}}

\nwc{\hhn}{H^{n-1}}
\nwc{\hhns}{H^{n-1}_{*}}
\nwc{\hhps}{H^{p}_{*}}
\nwc{\hhu}{H^{1}}
\nwc{\hhus}{H^{1}_{*}}
\nwc{\hhz}{H^{0}}
\nwc{\hhzs}{H^{0}_{*}}

\nwc{\nn}{\mathbb{N}}

\nwc{\oito}{\textbf{(viii)}}
\nwc{\oo}{\mathcal{O}}
\nwc{\oox}{\mathcal{O}_{X}}

\nwc{\pp}{\mathbb{P}}

\nwc{\qua}{\textbf{(iv)}}

\nwc{\seis}{\textbf{(vi)}}
\nwc{\sete}{\textbf{(vii)}}

\nwc{\tre}{\textbf{(iii)}}

\nwc{\um}{\textbf{(i)}}

\nwc{\ww}{\omega}
\nwc{\wwx}{\omega _{X}}

\nwc{\zz}{\mathbb{Z}}
\nwc{\C}{\mathbb{C}}

\newtheorem{coro}{Corollary}[section]
\newtheorem{defi}[coro]{Definition}

\newtheorem{lema}[coro]{Lemma}
\newtheorem{prop}[coro]{Proposition}

\newtheorem{teo}[coro]{Theorem}

\begin{document}

\title{Linear \& Steiner bundles \\ over projective varieties}

\author{Marcos Jardim}
\author{Renato Vidal Martins}

\address{IMECC - UNICAMP \\
Departamento de Matem\'atica \\ Caixa Postal 6065 \\
13083-970 Campinas-SP, Brazil}
\email{jardim@ime.unicamp.br}
\address{ICEx - UFMG \\
Departamento de Matem\'atica \\ 
Av. Ant\^onio Carlos 6627 \\
30123-970 Belo Horizonte MG, Brazil}
\email{renato@mat.ufmg.br}


\begin{abstract}
We use a generalization of Horrocks monads for arithmetic Cohen-Macaulay (ACM) varieties to
establish a cohomological characterization of linear and Steiner bundles over projective 
spaces and quadric hypersurfaces. We also study resolutions of bundles on ACM varieties by
line bundles, and characterize linear homological dimension in the case of quadric hypersurfaces.
\end{abstract}

\maketitle


\section{Introduction}

Monads have been introduced by G. Horrocks' \cite{H1} in the 60's. There he 
proved that every locally free sheaf $E$ on $\pp^n$ is the cohomology of a 
monad for which the first and last terms are sums of line bundles while the 
middle term, say $B$, satisfies $H^{1}(B(k))=H^{n-1}(B(k))=0$ for every 
$k\in\zz$ and also $H^{p}(B(k))\cong H^{p}(E(k))$ for every 
$2\leq p\leq n-2$ and $k\in\zz$. The idea of Horrocks 
was, in his own words, that one can always ``kill" cohomology of any 
locally free sheaf on a projective space. This motivated the study of some 
special kind of monads which were named \emph{Horrocks} later on (cf. 
\cite{AR} and \cite{BH} for example).

Monads on projective spaces have been much studied in the past 25 years (see for instance
\cite{F,J-i,MPR,OSS} and the references therein). More recently, many authors 
have also been interested on monads over more general varieties, see 
\cite{B,CMR,JMR}. We start the present article in Section \ref{secmon} studying general monads in
arbitrary varieties, providing necessary and sufficient conditions that a coherent
sheaf must satisfy in order to be the cohomology of a monad (cf. Proposition \ref{prpnec}). 
When so, it is easy to see that there will be infinitely many monads for each coherent
sheaf, so we introduce a notion of equivalence between monads.

In Section \ref{sechor}, we prove a generalization to ACM varieties of the above mentioned result
due to Horrocks. Although it is not available in the literature, such generalization might be known
to specialists. It will be our main tool in the cohomological characterization of linear and Steiner
bundles.

We devote Section \ref{seclin} to the study of linear bundles. Linear bundles on 
$\pp^2$ and $\pp^3$ have been studied since the late 1970's (cf. \cite{OSS}), and the 
mathematical instanton bundles on $\pp^{2n+1}$ introduced by Okonek and Spindler in
\cite{OS} are examples of linear bundles. More general linear monads were first considered
in \cite{JMR}. A cohomological characterization of linear torsion free sheaves on $\pp^n$
has been recently obtained in \cite{J-i} and L. Costa and R. M. Mir\'o-Roig have also obtained
a cohomological characterization of linear bundles of rank $2n$ on quadric 
hypersurfaces on $\pp^{2n+2}$ ($n\ge2$) (cf. \cite[Theorem 4.5]{CMR}). Both 
results were obtained by means of Beilinson spectral sequence. Here we prove 
a similar result for projective spaces (cf. Theorem \ref{thmprj}) and a 
generalization for the case of quadric hypersurfaces (cf. Theorem \ref{thmhyp}) 
with more elementary tools, essentially just using the generalized Horrocks' theorem and the
splitting criteria for locally free sheaves on these varieties.

In Section \ref{secstn} we study Steiner bundles. They were introduced by Dolgachev 
and Kapranov in the context of hyperplane arrangements on $\pp ^n$ (cf. \cite{DK}), 
and have been studied by several other authors since then, see for instance \cite{AO,BS,Br}
for Steiner bundles on projective spaces and \cite{S} for Steiner bundles on quadrics.
Here we apply the same technique we used to characterize linear 
bundles, once again without using Beilinson spectral sequence, in order to obtain 
a characterization of Steiner bundles on $\pp ^{n}$ (cf. Theorem \ref{thmstp}) 
which is slightly different from \cite[Proposition 3.2]{DK}. Besides, after naturally 
extending the concept of Steiner bundles to arbitrary varieties, we also 
give a characterization of such bundles on quadric hypersurfaces (cf. Theorem \ref{thmstq}).

We finish this work in Section \ref{secres} studying linear homological dimension (lhd) 
of vector bundles, i.e., the minimal length of a locally-free resolution whose terms split
as sum of line bundles, a notion introduced by G. Bohnhorst and H. Spindler in \cite{BS}.
They proved (cf. \cite[Proposition 1.4]{BS}) that if $E$ is a vector bundle on $\pp^n$ 
satisfying $H^p_*(E)=0$ for $1\le p\le n-d-1$, then ${\rm lhd}(E)\le d$. Here we prove a
generalization of this result for general ACM varieties (cf Proposition \ref{longacm}) and a 
version for on quadric hypersurfaces (cf. Theorem \ref{thmres}): a locally free sheaf $E$
on $Q_n$ ($n\ge 3$) has ${\rm lhd}(E)\le d$ if and only if $H^p_*(E)=0$ for all
$1\le p\le n-d-1$ and $H^{n-d-1}_*(E\otimes\Sigma)=0$ where $\Sigma$ is any spinor bundle
on $Q_{n}$.
  
\paragraph{\bf Acknowledgments.}
The first named author is partially supported by FAEPEX grant number 1652/04 
and the CNPq grant number 300991/2004-5. The second author is partially supported by CNPq grant number PDE 200999/2005-2. We thank N. Mohan Kumar and 
Rosa Maria Mir\'o-Roig for answering of our questions regarding their work.


\section{Monads on Projective Varieties}
\label{secmon}

Let $X$ be a \emph{projective variety}, i.e., a projective scheme over an 
algebraically closed field together with a given very ample invertible sheaf 
which we denote by $\oox(1)$. For any coherent sheaf $E$ on $X$ we denote 
$E(k)=E\otimes\oox(k)$ and $H^{p}_{*}(E):=\oplus _{k\in\zz}\ H^{p}(E(k))$. 
In particular, the graded ring of homogeneous functions on $X$ is given by 
$S(X)=H^{0}_{*}(\oox)$. Note that $H^{p}_{*}(E)$ is a graded $S(X)$-module. 
Finally, $\wwx$ denotes the dualizing sheaf on $X$.  

\begin{defi}
\label{defmon}
A \emph{monad} on a projective variety $X$ is a complex 
$$ 
{\rm M}_\bullet ~~:~~ 0 \longrightarrow A \stackrel{\alpha}{\longrightarrow} 
B \stackrel{\beta}{\longrightarrow} C\longrightarrow 0
$$
of locally free sheaves $A$, $B$ and $C$ on $X$ which is exact on the first 
and last terms. The coherent sheaf $E=\ker\beta/{\rm Im}\alpha$ is called the 
\emph{cohomology} of ${\rm M}_\bullet$ and one also says that ${\rm M}_\bullet$ 
is a \emph{monad for} $E$. 
\end{defi}

Every monad on $X$ can be broken down, using the fact that $\alpha$ is 
injective and $\beta$ is surjective, into two short exact sequences
\begin{equation}
\label{equker}
0\longrightarrow K\longrightarrow B\stackrel{\beta}{\longrightarrow} C 
\longrightarrow 0\ \ \text{and} 
\end{equation}
\begin{equation}
\label{equimg}
0\longrightarrow A\stackrel{\alpha}{\longrightarrow} K\longrightarrow E
\longrightarrow 0
\end{equation}
where $K=\ker\beta$ is also locally-free. 

To a given monad ${\rm M}_\bullet : 0 \rightarrow A \stackrel{\alpha}{\rightarrow} 
B \stackrel{\beta}{\rightarrow} C\rightarrow 0$ whose cohomology is $E$ we can 
also associate the \emph{dual} monad 
${\rm M}_\bullet ^* : 0 \rightarrow C^* \stackrel{\beta ^*}{\rightarrow} 
B^* \stackrel{\alpha ^*}{\rightarrow} A^*\rightarrow 0$ whose cohomology is 
precisely $E^*$. It is easy to see that not every coherent sheaf on $X$ can be 
obtained as the cohomology of a monad. We have the following conditions:

\begin{prop}
\label{prpnec}
If $E$ is a coherent sheaf on a projective variety $X$, then:
\begin{itemize}
\item[\um] there exists a monad for $E$ if and only if 
${\mathcal E}xt^p(E,F)=0$ for every sheaf of $\oox$-modules $F$  and each 
$p\geq2$;
\item[\doi] if $X$ is smooth then there exists a monad for $E$ if and only if 
${\mathcal E}xt^p(E,\oox )=0$ for each $p\geq2$; 
\item[\tre] if $X$ is Cohen-Macaulay with pure dimension $n\geq 3$ then the 
existence of a monad for $E$ implies that for each $1\leq p\leq n-2$, 
$H^p_*(E)$ is finitely generated as a graded $S(X)$-module.
\end{itemize}
\end{prop}

\begin{proof}
If there is a monad for $E$, we can break it down like above to obtain the 
exact sequence (\ref{equimg}) where $A$ and $K$ are both locally free. 
Conversely, if we find an exact sequence like (\ref{equimg}) for $E$, then $E$ 
is the cohomology of the following monad:
$$ 
0 \longrightarrow A \longrightarrow K \longrightarrow 0 \longrightarrow 0. 
$$
In other words, there exists a monad for $E$ if and only if $E$ admits a 
locally free resolution of length at most $1$. Therefore (i) and (ii) follows 
from \cite[pp. 127-128]{M}.

To prove (iii), regarding (\ref{equimg}) again, we note that if $X$ fulfills 
the requirements above, then, by duality, there are integers $m$, $M$ such 
that for each $1\leq p\leq n-1$, we have $H^p(A(k))=H^p(K(k))=0$ for all 
$k<m$ and all $k>M$ since $A$ and $K$ are locally free. From the cohomology 
sequence associated to (\ref{equimg}) 
$$ 
H^p(K(k)) \to H^p(E(k)) \to H^{p+1}(A(k)) 
$$
we are able to conclude that, for $1\leq p\leq n-2$, $H^p(E(k))=0$ for all 
$k<m$ and $H^p(E(k))=0$ for all $k>M$. It follows that $H^p_*(E)$ is finitely 
generated for $1\leq p\leq n-2$.
\end{proof}

It is also easy to see that it can exist many monads for a given coherent 
sheaf which satisfies the first item of the above result. So it is important 
to introduce a notion of equivalence between monads. 

\begin{defi}
We say that two monads are \emph{equivalent} if their associated exact 
sequences (\ref{equker}) and (\ref{equimg}) are isomorphic as extensions.
\end{defi}

We finish this section just noticing that the numerical invariants of the 
cohomology sheaf $E$ of a monad $0\to A\to B\to C\to 0$ are easily computable 
in terms of the numerical invariants of the sheaves $A$, $B$ and $C$. More
precisely, the Chern character of $E$ is given by
\begin{equation}
\label{equchn}
ch(E)=ch(B)-ch(A)-ch(C).
\end{equation}
In particular, the rank of $E$ is given by
\begin{equation}
\label{equrnk}
rk(E)=rk(B)-rk(A)-rk(C)
\end{equation}
which is another relation which will also be useful in the next section.


\section{Horrocks Monads on ACM Varieties}
\label{sechor}

>From now on we restrict our studies of arbitrary monads, varieties and sheaves 
to, respectively, Horrocks, ACM and locally free ones. More precisely, we start
with the following definition (cf. \cite[Definition 2.2]{AR} and \cite[Section 3]{BH}): 

\begin{defi}
\label{defhor}
A monad 
\begin{equation}\label{mon.def}
{\rm M}_{\bullet}\ :\ 0 \to A \stackrel{\alpha}{\longrightarrow} B 
\stackrel{\beta}{\longrightarrow} C \to 0
\end{equation}
on a projective variety $X$ whose cohomology is $E$ is called \emph{Horrocks} if
\begin{itemize}
\item[\um ] $A$ and $C$ are direct sums of invertible sheaves;
\item[\doi ] $\hhus (B)=\hhns (B)=0$ and if $n\geq 4$ then $H^{p}_{*}(B)\cong H^{p}_{*}(E)$ for $2\leq p\leq n-2$;
\end{itemize}
if besides the monad satisfies
\begin{itemize}
\item[\tre ] no direct summand of $A$ is isomorphic to a direct summand of $B$;
\item[\qua ] no direct summand of $C$ is the image of a line subbundle of $B$;
\end{itemize}
then it is also called \emph{minimal}.
\end{defi}

It is relevant to ask whether we have uniqueness for minimal monads 
for a given sheaf. Later on in this section we will guarantee 
this property for some kind of monads. Now we specify the varieties we will be interested in. 

\begin{defi} 
A projective variety $X$ of pure dimension $n$ is \emph{arithmetically Cohen-Macaulay} 
(ACM) if $H^{p}_{*}(\oox )=0$ for every $1\leq p\leq n-1$.
\end{defi}

This definition is equivalent to saying that if $\oo _{X}(1)$ corresponds to an
embedding $X\subset\pp ^{r}$ then $H^{1}_{*}(\pp^r,{\mathcal I}_X )=0$ and that
its homogeneous coordinate ring $S(X)$ is Cohen-Macaulay \cite{CH}. Every complete 
intersection scheme $X\subset\pp^r$ is ACM. Nonsingular Fano varieties with 
cyclic Picard group are also ACM.

\begin{teo} 
\label{thmhor}
Let $X$ be an ACM variety of dimension $n\ge3$ and let $E$ be a locally free 
sheaf on $X$. Then there is a 1-1 correspondence between collections 
$$
\{h_{1},\ldots ,h_{r},g_{1}, \ldots ,g_{s}\}\ \text{with}\ 
h_{i}\in H^{1}(E^*\otimes\wwx (k_{i}))\ \text{and}\ g_{j}\in\hhu (E(-l_{j}))
$$ 
for integers $k_{i}$'s and $l_{j}$'s and
equivalence classes of monads for $E$ of the form
$$ 
{\rm M}_{\bullet}\ :\ 0\longrightarrow\ \oplus _{i=1}^{r}\ \wwx (k_{i})\ 
\stackrel{\alpha}{\longrightarrow}\ F\ \stackrel{\beta}{\longrightarrow}\ 
\oplus _{j=1}^{s}\ \oox (l_{j})\ \longrightarrow 0.
$$
This correspondence is such that:
\begin{itemize}
\item[\um ]  ${\rm M}_{\bullet}$ is Horrocks if and only 
if the $g_{j}$'s generate $\hhus (E)$ and the $h_{i}$'s generate 
$H^{1}_{*}(E^*\otimes\wwx )$ as $S(X)$-modules;
\item[\doi ] ${\rm M}_{\bullet}$ is minimal Horrocks if and only if the 
$g_{j}$'s constitute a minimal set of generators for $\hhus (E)$ and the 
$h_{i}$'s constitute a minimal set of generators for $H^1_*(E^*\otimes\wwx)$ 
as $S(X)$-modules.
\end{itemize}
\end{teo} 

\begin{proof}
The element 
$\sum _{i} h_{i}\in\oplus _{i}\ H^{1}(E^{*}\otimes\wwx (k_{i}))\cong
\text{Ext}^{1}(E,\oplus _{i}\ \wwx(k_{i}))$ yields an extension
\begin{equation}
\label{equex1}
0\longrightarrow\ \oplus _{i=1}^{r}\ \wwx (k_{i})\ 
\stackrel{\alpha}{\longrightarrow}\ K\ \longrightarrow\ E\ \longrightarrow 0.
\end{equation}

On the other hand, since $X$ is ACM, we have that  
$H^{p}(K(m))\cong H^{p}(E(m))$ for
$m\in\zz$ and $1\leq p\leq n-2$. Set $g_{j}'\in H^{1}(K(-l_{j}))$ to be the 
image of $g_{j}$ for each $j$. Then 
$\sum _{j} g_{j}'\in\oplus _{j}\ H^{1}(K(-l_{j}))\cong
\text{Ext}^{1}(\oplus _{j}\ \oox(l_{j}),K)$ yields an extension
\begin{equation}
\label{equex2}
0\longrightarrow\ K\ \longrightarrow\ F\ \stackrel{\beta}{\longrightarrow}\ 
\oplus _{j=1}^{s}\ \oox(l_{j})\ \longrightarrow 0.
\end{equation}
Now (\ref{equex1}) and (\ref{equex2}) yield the monad  
\begin{equation}
\label{equmon}
0\longrightarrow\ \oplus _{i=1}^{r}\ \wwx (k_{i})\ 
\stackrel{\alpha}{\longrightarrow}\ F\ \stackrel{\beta}{\longrightarrow}\ 
\oplus _{j=1}^{s}\ \oox (l_{j})\ \longrightarrow 0
\end{equation}
of which $E$ is the cohomology.

Conversely, given (an equivalence class of) a monad like (\ref{equmon}), 
we can break it down into two exact sequences as in (\ref{equex1}) and 
(\ref{equex2}) where $K=Ker \beta$. Then (\ref{equex1}) corresponds to an 
element $h\in\text{Ext}^{1}(E,\oplus _{i}\ \wwx(k_{i}))\cong\oplus _{i} 
H^{1}(E^*\otimes\wwx (k_{i}))$ and we write $h=\sum _{i} h_{i}$ in order 
to have the desired $h_{i}$'s. The same can be done with (\ref{equex2}) to get 
back the $g_{j}$'s.

Now, since $\hhus (\oplus _{j}\ \oox (l_{j}))=0$, 
from (\ref{equex2}) we have that $\hhus (F)=0$ if and only if  for every 
$m\in\zz$ the boundary morphisms 
$\psi _{m}: \hhz (\oplus _{j}\oox (l_{j}+m))\longrightarrow
\hhu (K(m))$ are surjective. But now we have 
\begin{gather*}
\begin{matrix}
\psi _{m}& : &\oplus _{j=1}^{s}\hhz (\oox (l_{j}+m)) & \longrightarrow & 
H^{1}(K(m))\\
         &   &\sum _{j} f_{j}& \longmapsto &\sum _{j} f_{j}g_{j}'
\end{matrix}
\end{gather*}
where the $f_{j}$'s are polynomials of degree $l_{j}+m$. Then all the 
$\psi _{m}$'s are surjective if and only if the $g_{j}'$'s generate 
$\hhus (K)$, what happens if and only if the $g_{j}$'s generate $\hhus (E)$
as $S(X)$-module.

We also have that $H^{n-1}_*(E)$ vanishes if and only if 
$H^{1}_*(E^*\otimes\wwx )$ does. Tensoring the dual monad 
${\rm M}_{\bullet}^{*}$ by $\wwx$ we conclude from the same argument above 
that the latter assertion holds if and only if the $h_{i}$'s generate 
$H^1_*(E^{*}\otimes\wwx )$ as $S(X)$-module.

Since we have  
$H^{p}(F(m))\cong H^{p}(K(m))\cong H^{p}(E(m))$ for every $2\leq p\leq n-2$ and
$m\in\zz$ we are done for (i).

To prove (ii), for each $l$ let us gather the $g_j$'s which are elements of $H^1(E(-l))$, say 
they are $g_1^l,\ldots ,g_{s_l}^l$. If we identify $H^1_*(K)=H^1_*(E)$ by means of 
(\ref{equex1}) then from (\ref{equex2}) we have for every $l$ a morphism
\begin{gather*}
\begin{matrix}
\psi _{l}& : &\oplus _{j=1}^{s_l}\hhz (\oox ) & \longrightarrow & 
H^{1}(E(-l))\\
         &   &(c_1,\ldots ,c_{s_l})& \longmapsto &\sum _{j=1}^{s_l}c_{j}g_{j}^l .
\end{matrix}
\end{gather*}
We have that the $g_j$'s constitute a minimal set of generators for $H^1_*(E)$ if and only if the $\psi _l$'s are
isomorphisms. Since the $\psi _l$'s were induced from (\ref{equex2}) this is equivalent to the fact that no 
$\oox (l_j)$ is the image of a line subbundle of $F$. In order to conclude something similar for the $h_i$'s we just 
consider the dual monad and apply the same argument.
\end{proof}

In what follows we will derive some interesting consequences of the above result.
We start by generalizing Horrocks' theorem to ACM varieties.

\begin{coro}
\label{corhor}
Every locally free sheaf $E$ on an ACM variety of dimension $n\geq 3$ is the 
cohomology of a monad
$$ 0\to A\to B\to C\to 0 $$
where $A$ and $C$ split as sums of line bundles and $B$ satisfies:
\begin{itemize}
\item[\um] $H^{1}_{*}(B)=H^{n-1}_{*}(B)=0$;
\item[\doi] for $n\geq 4$, $H^{p}_{*}(E(k))\cong H^{p}_{*}(B(k))$ for
$2\leq p\leq n-2$.
\end{itemize} 
\end{coro}

The proof is immediate from Theorem \ref{thmhor} and the fact that both 
$H^{1}_{*}(E)$ and $H^{1}_{*}(E^*\otimes\wwx )$ are finitely generated as 
$S(X)$-modules because $E$ is locally free. The corollary supports the
conclusion that {\em line bundles and locally-free sheaves 
$B$ satisfying $H^1_*(B)=0=H^{n-1}_*(B)$ are building blocks for more 
general locally-free sheaves on ACM varieties}.

Now consider the following definition.

\begin{defi}
A coherent sheaf $W$ on an ACM variety of dimension $n\ge2$ is said to be 
{\em arithmetically Cohen-Macaulay} (ACM) if $H^p_*(W)=0$ for $1\le p\le n-1$ 
and if there exist integers $s$ and $t$ such that $H^0(W(k))=0$ for $k\leq s$ 
and $H^{n}(W(k))=0$ for $k\geq t$. Taking $s$ (resp. $t$) to be the greatest
(resp. smallest) integer satisfying the related property, we say that $W$ is an
ACM sheaf \emph{with parameters} $s$ \emph{and} $t$.
\end{defi}

As it was observed in \cite[Proposition 2.1]{CH}, ACM sheaves correspond to 
graded maximal Cohen-Macaulay modules on $S(X)$. We therefore obtain the 
following statement:

\begin{coro}
Every locally-free sheaf on an ACM variety of dimension $3$ 
is the cohomology of a monad
$$ 
0 \to A \to W \to C \to 0 
$$
where $A$ and $C$ are sums of line bundles and $W$ is an ACM sheaf.
\end{coro}

The natural question then becomes how to classify indecomposable locally free 
ACM sheaves on a given variety, which can be a very hard problem.

As it is well known, the only indecomposable locally free ACM sheaves on 
$\pp^n$ are line bundles. On quadric hypersurfaces $Q_n\subset\mathbb{P}^n+1$ 
($n\ge3$), according to \cite{K}, the only indecomposable locally-free ACM 
sheaves are line bundles and twists of the spinor bundles by line 
bundles. In particular, we can conclude that:

\begin{coro}
Every locally-free sheaf on a quadric hypersurface $Q_3\subset\mathbb{P}^4$
is the cohomology of a monad
$$ 
0 \to A \to W \to C \to 0 
$$
where $A$ and $C$ are sums of line bundles and $W$ is a sum of line bundles and
twisted spinor bundles.
\end{coro}

On the other hand, hyperplanes and quadrics are the only hypersurfaces in 
projective space for which there are only a finite number of indecomposable 
locally-free ACM sheaves, up to twisting by a line bundle (cf. \cite{BGS}). 
Some specific varieties have been looked at in the literature; in \cite{AG} 
the authors classify all locally-free ACM sheaves on the grassmannian of lines 
in $\mathbb{P}^4$; certain Fano 3-folds were considered in \cite{AC}.

\

Theorem \ref{thmhor} also lead us naturally to analyze the special case of 
monads which appear in the following result.

\begin{prop}
\label{prplmm}
Let $X$ be an ACM variety of dimension $n\ge2$. If $E$ is the cohomology of a
monad of the form
$$ 
{\rm M}_{\bullet}\ :\ 0\longrightarrow\ \oox (l)^{\oplus a}\ 
\longrightarrow\ \oox ^{\oplus b}\ \longrightarrow\ 
\oox (m)^{\oplus c}\longrightarrow 0
$$
with $l<0<m$ then:
\begin{itemize}
\item[\um] ${\rm ch}(E)=b-a\cdot{\rm ch}(\oox(l))-c\cdot{\rm ch}(\oox(m))$;
\item[\doi] $b=rk(E)-a-c$ and if $\wwx =\oox(\lambda)$ then
$a=h^{n-1}(E(\lambda -l))=h^1(E^*(l))$ and $c=h^{1}(E(-m))$;
\item[\tre] ${\rm M}_{\bullet}$ is, up to isomorphism, the unique monad for $E$
of the stated form. In particular, $l$ and $m$ are also determined by $E$;
\item[\qua ] ${\rm M}_{\bullet}$ is minimal Horrocks. 
\end{itemize}
if $W$ is an ACM sheaf on $X$ with parameters $s$ and $t$ then:
\begin{itemize}
\item[\cin ] if $n\ge2$, $H^0(E\otimes W(s))=H^n(E\otimes W(t))=0$;
\item[\seis ] if $n\ge3$, $H^1(E\otimes W(s-m))=H^{n-1}(E\otimes W(t-l))=0$;
\item[\sete ] if $n\ge4$, $H^p_*(E\otimes W)=0$ where $2\le p\le n-2$.
\end{itemize}
\end{prop}

\begin{proof}
We have that (i) follows from (\ref{equchn}) and the first assertion of (ii)
follows from (\ref{equrnk}). Since $X$ is ACM we can adapt 
the proof of \cite[Proposition 4]{BH} with the remark which appears right 
after to conclude (iii) taking into account that 
$Hom(\oox ,\oox (l))=Hom(\oox (m),\oox )=0$ if $l<0<m$. In order to prove
(v), (vi) and (vii), consider the sequences (\ref{equker}) and
(\ref{equimg}) twisted by $W(k)$:
$$ 
0 \to K\otimes W(k) \to W(k)^{\oplus b} \to W(k+m)^{\oplus c} \to 0\ \ \ 
{\rm and}
$$
$$ 
0 \to W(k+l)^{\oplus a} \to  K\otimes W(k) \to E\otimes W(k) \to 0.
$$
Passing to cohomology, one obtains the desired results. From (vii) with 
$W=\oox$ we see that ${\rm M}_{\bullet}$ is clearly Horrocks and from the very
definition of minimality it is also minimal and so (iv) is proved. Finally, 
the second assertion of item (ii) is a consequence of (iv) and Theorem \ref{thmhor}.
\end{proof}


\section{Linear Bundles}
\label{seclin}

In this section we consider a particular case of the monads appearing in 
Proposition \ref{prplmm} the ones with $l=-1$, $m=1$.

\begin{defi}
\label{deflin}
A monad on a projective variety $X$ is called \emph{linear} if it is of the 
following form: 
\begin{equation}
\label{equlin}
0 \to \oox(-1)^{\oplus a} \to \oox^{\oplus b} \to \oox(1)^{\oplus c} \to 0.
\end{equation}
Similarly the cohomology of a linear monad is called a \emph{linear} sheaf.
\end{defi}

Our goal is to show that the converse of Proposition \ref{prplmm} holds
in the case of linear monads over projective spaces and quadric hypersurfaces,
characterizing linear locally-free sheaves over these varieties. One of our main
tools will be the following lemma.

\begin{lema} 
\label{lemgen}
Let $X$ be a non-singular ACM variety of dimension $n\ge2$. If $E$ is a locally
free sheaf on $X$ for which there exists an $l$ such that $H^{p}(E(l-p+1))=0$ 
for $2\leq p\leq n$ then $H^{1}(E(l))$ generates 
$\oplus _{k\geq l}\ H^{1}(E(k))$ as an $S(X)$-module.
\end{lema}

\begin{proof}
Let $H$ be a hyperplane divisor on $X$. Using the restriction sequence
\begin{equation}
\label{equhyp}
0 \to E(k-1) \to E(k) \to E(k)|_H \to 0 
\end{equation}
we can broaden the hypothesis above into: 
\begin{equation}
\label{equbro}
H^{p}(E(k))=0\ \text{for}\ k\geq l-p+1\ \text{and}\ 2\leq p\leq n
\end{equation}

Since $X$ is smooth, we can apply Bertini's Theorem successively to obtain
a regular sequence $\{x_{1},\ldots ,x_{n}\}$ for $S(X)$. Now set 
$Y:=Z(x_{1},\ldots ,x_{n})$ to be the common zero locus. The Koszul complex of 
the above sequence yields the following chain:
$$ 
0 \to \oo _{X}(-n) \to \wedge ^{n-1}\oo_{X}^{\oplus n}(-n+1) \to\ldots\to
\wedge ^{2}\oo_{X}^{\oplus n}(-2) 
\to\oo_{X}^{\oplus n}(-1)\to\mathcal{I}_{Y}\to 0. 
$$

Twisting it by $E(k)$, using exterior algebra properties, we are led to
$$ 
0 \to E(k-n)\stackrel{\alpha _{n}^{k}}{\longrightarrow}
E(k-n+1)^{\oplus d_{n-1}}\stackrel{\alpha _{n-1}^{k}}{\longrightarrow}\ldots
\stackrel{\alpha _{2}^{k}}{\longrightarrow} 
E(k-1)^{\oplus n}\stackrel{\alpha _{1}^{k}}{\longrightarrow}
E\otimes\mathcal{I}_{Y}(k)\to 0. 
$$
where $d_{p}={n\choose p}$. So we have the following short exact sequences
\begin{gather*}
\begin{matrix}
  & 0 \longrightarrow  Ker\ \alpha _{1}^{k}  \longrightarrow  
E(k-1)^{\oplus n}\longrightarrow  E\otimes\mathcal{I}_{Y}(k)  
\longrightarrow  0 & \\
  & \ldots         & \\
  & 0 \longrightarrow  Ker\ \alpha _{p}^{k}  \longrightarrow  
E(k-p)^{\oplus d_{p}}\longrightarrow  Ker\ \alpha _{p-1}^{k}  \longrightarrow  
0 & \\
  & \ldots         & \\
  & 0 \longrightarrow  E(k-n)  \longrightarrow  
E(k-n+1)^{\oplus d_{n-1}} \longrightarrow  Ker\ \alpha_{n-2}^{k}  
\longrightarrow  0.& 
\end{matrix}
\end{gather*}

Considering the long exact sequences for cohomology and using (\ref{equbro}) 
we have that
\begin{itemize}
\item[ ] $H^p(Ker\ \alpha _{p-1}^{k})\cong H^{p+1}(Ker\ \alpha _{p}^{k})$\ \ \
for $k\geq l+1$ and $2\leq p\leq n-2$,
\item[ ] $H^{n-1}(Ker\ \alpha _{n-2}^{k})=0$\ \ \ for $k\geq l+1$,
\end{itemize}
and so we can conclude that $H^2(Ker\ \alpha _{1}^{k})=0$ for $k\geq l+1$. 
In particular 
\begin{equation}
\label{equsj1}
H^{1}(E(k-1))^{\oplus n}\longrightarrow H^{1}(E\otimes\mathcal{I}_{Y}(k))
\longrightarrow 0
\end{equation}
is a surjection for every $k\geq l+1$. Now if we consider the exact sequence 
$$ 
0\longrightarrow E\otimes\mathcal{I}_{Y}(k)\longrightarrow E(k)\longrightarrow
E\otimes\oo _{Y}(k)\longrightarrow 0 
$$
we also have a surjection
\begin{equation}
\label{equsj2}
H^{1}(E\otimes\mathcal{I}_{Y}(k))\longrightarrow H^{1}(E(k))\longrightarrow 0
\end{equation}
because $Y$ is zero dimensional. Therefore, combining (\ref{equsj1}) with 
(\ref{equsj2}), we have a surjection
$$ 
H^{1}(E(k-1))^{\oplus n}\longrightarrow H^{1}(E(k))\longrightarrow 0 
$$
for every $k\geq l+1$. So if $\{g_{1},\ldots ,g_{s}\}$ is a basis of the 
vector space $H^{1}(E(l))$ then, by construction, every element 
$h\in H^{1}(E(l+1))$ can be written as 
$h=\sum _{j=1}^{s}(\sum _{i=1}^{n}c_{ij}x_{i})\ g_{j}$ where the $c_{ij}$'s are
constants. Reapplying the argument successively for every $k\geq l+1$ we see 
that $\{g_{j}\}_{j=1}^{s}$ generate $\oplus _{k\geq l}\ H^{1}(E(l))$ as an 
$S(X)$-module. 
\end{proof}

\subsection{Linear Bundles on Projective Spaces}

Gathering Theorem \ref{thmhor}, Proposition \ref{prplmm} and Lemma \ref{lemgen}
we obtain an alternative, more elementary proof of a result proved in \cite{J-i}
using the Beilinson spectral sequence:

\begin{teo}
\label{thmprj}
Let $E$ be a locally free sheaf on $\pp^n$ where $n\geq 3$. Then $E$ is a linear
sheaf realized as the cohomology of a monad
$$ 
0 \to \oo_{\pp^n}(-1)^{\oplus a} \to \oo_{\pp^n}^{\oplus b}
\to \oo_{\pp^n}(1)^{\oplus c} \to 0 
$$
if and only if the following cohomological conditions hold:
\begin{itemize}
\item[\um] ${\rm ch}(E)=b-a\cdot{\rm ch}(\oo _{\pp ^{n}}(-1))
-c\cdot{\rm ch}(\oo _{\pp ^{n}}(1))$;
\item[\doi] $H^0(E(-1))=H^1(E(-2))=H^{n-1}(E(1-n))=H^n(E(-n))=0$;
\item[\tre] for $n\ge4$, $H^p_*(E)=0$ where $2\le p\le n-2$.
\end{itemize}
\end{teo}

\begin{proof}
The ``only if" part is clear from Proposition \ref{prplmm} taking 
$W=\oo_{\pp^n}$. Conversely since $H^n(E(-n))=H^{n-1}(E(1-n))=0$ and
$H^p_*(E)=0$ for $2\leq p\leq n-2$ we conclude from Lemma \ref{lemgen} that
$H^{1}(E(-1))$ generates $\oplus _{k\geq -1}\ H^{1}(E(k))$ as an 
$S(\pp^n)$-module. Now since $H^{0}(E(-1))=H^{1}(E(-2))=0$ we have from 
(\ref{equhyp}) that $H^{1}(E(k))=0$ for $k\leq -2$. Thus $H^{1}(E(-1))$
actually generates $H^{1}_{*}(E)$ as $S(\pp^n)$-module. Similarly by duality 
$E^{*}$ satisfies the same conditions (ii) and (iii) above, hence $H^{1}(E^{*}(-1))$
generates $H^{1}_{*}(E^{*})$ as an $S(\pp^n)$-module. Conditions (i)-(iii) also 
imply that $c=-\chi (E(-1))=h^{1}(E(-1))$ and 
$a=-\chi (E^*(-1))=h^{1}(E^*(-1))$.

Therefore from Theorem \ref{thmhor} we have that $E$ is the cohomology of a 
monad of the form
$$ 
0\longrightarrow\oo _{\pp ^n}(-1)^{\oplus a}\longrightarrow F
\longrightarrow\oo _{\pp ^{n}}(1)^{\oplus c}\longrightarrow 0.
$$
with $H^p_*(F)=0$ for $1\leq p\leq n-1$. Thus by Horrocks' splitting criterion 
$F$ must split as a sum of line bundles. Moreover, (i) implies that $rk(F)=b$
and ${\rm ch}_{2}(F)=0$. So if we set $F=\oplus_{i=1}^{b} \oo_{\pp^n}(a_i)$ it 
follows that ${\rm ch}_2(F)=\frac{1}{2}\sum_{i=1}^{b} a_i^2=0$ and hence the 
$a_i's$ must vanish what implies $F=\oo_{\pp^n}^{\oplus b}$ as desired.
\end{proof}

Linear bundles on $\pp^n$ exist if and only if at least either $b\geq 2c+n-1$ and 
$b\geq a+c$ or $b\geq a+c+n$ (cf. \cite[Main Theorem]{J-i}).

\subsection{Linear bundles on quadric hypersurfaces}

Next we obtain a characterization of linear sheaves on quadric hypersurfaces 
$Q_{n}$ on $\pp^{n+1}$ ($n\ge4$), generalizing a result obtained in 
\cite[Theorem 4.5]{CMR} using the Beilinson spectral sequence. We will use
some basic facts regarding spinor bundles on quadrics, see \cite{Ot,O1}. 

\begin{teo}\label{thmhyp}
Let $E$ be a locally free sheaf on a quadric hypersurface $Q_{n}$ on 
$\pp ^{n+1}$ where $n\ge3$. Then $E$ is a linear sheaf realized as
the cohomology of a monad
$$ 
0\to\oo_{Q_{n}}(-1)^{\oplus a}\to\oo_{Q_{n}}^{\oplus b}
\to\oo_{Q_{n}}(1)^{\oplus c}\to 0 
$$
if and only if the following cohomological conditions hold:
\begin{itemize}
\item[\um ] ${\rm ch}(E)=b-a\cdot{\rm ch}(\oo _{Q_{n}}(-1)) -c\cdot{\rm ch}(\oo _{Q_{n}}(1))$;
\item[\doi ] $H^0(E(-1))=H^1(E(-2))=H^{n-1}(E(2-n))=H^n(E(1-n))=0$;
\item[\tre ] if $n\geq 4$, $H^p_*(E)=0$ for each $2\le p\le n-2$;
\item[\qua ] if $n=3$,  $H^0(E\otimes\Sigma)=H^1(E\otimes\Sigma (-1))=H^2(E\otimes\Sigma (-1))=H^3(E\otimes\Sigma(-2))=0$
where $\Sigma$ is the spinor bundle on $Q_{3}$;
\item[\cin ] if $n\geq 4$, $H^p_*(E\otimes\Sigma)=0$ for each $2\leq p\leq n-2$ and every spinor bundle $\Sigma$ on $Q_n$. 
\end{itemize}
\end{teo}

\begin{proof} From Proposition \ref{prplmm} we have that if $E$ is the 
cohomology of the stated monad then (i)-(v) hold because $\oo _{Q_n}$ is an ACM sheaf 
with parameters $s=-1$ and $t=1-n$ and the spinor bundles
$\Sigma$ are ACM sheaves with parameters $s=0$
and $t=1-n$ (cf. \cite[Theorem 2.8]{O1}). 

Conversely, note that (ii)-(iii) together with Lemma \ref{lemgen} with $l=0$
imply that $H^1(E)$ generates $\oplus _{k\geq 0}\ H^{1}(E(k))$ as an 
$S(Q_n)$-module. Now since $H^{0}(E(-1))=H^{1}(E(-2))=0$ we have from 
(\ref{equhyp}) that $H^{1}(E(k))=0$ for $k\leq -2$. Thus $H^{1}(E(-1))\oplus H^{1}(E)$
generates $H^{1}_{*}(E)$ as $S(Q_n)$-module. Similarly by duality 
$E^{*}$ satisfies the same conditions (ii)-(iii) above, hence
$H^{1}(E^{*}(-1))\oplus H^{1}(E^*)$ generates $H^{1}_{*}(E^{*})$ as an $S(X)$-module.  

Therefore from Theorem \ref{thmhor} we have that $E$ is the cohomology of a 
monad of the form
\begin{equation}
\label{equofo} 
0\longrightarrow\oo _{Q_n}(-1)^{\oplus a}\oplus\oo _{Q_n}^{\oplus a'}\longrightarrow F
\longrightarrow\oo _{Q_n}^{\oplus c'}\oplus\oo _{Q_{n}}(1)^{\oplus c}\longrightarrow 0.
\end{equation}
with $H^p_*(F)=0$ for $1\leq p\leq n-1$. So $F$ is an ACM sheaf on $Q_{n}$. 

Breaking down the monad above we obtain two short exact
sequences as (\ref{equker}) and (\ref{equimg}). Tensoring them by a 
spinor bundle and using (v), we conclude that if $n\geq 4$ then $H^p_*(F\otimes\Sigma)=0$ for 
$2\le p\le n-2$ and every spinor bundle $\Sigma$. 

Now for $n=2m+1$ ($m\ge1$) tensoring the exact sequence of \cite[Theorem 2.8(i)]{O1}
by $F(k)$ we have:
\begin{equation}\label{fsigma}
0 \to F\otimes\Sigma (k)\to F(k)^{\oplus 2^{m+1}} \to F\otimes\Sigma(k+1) \to 0
\end{equation}
where $\Sigma$ is the only spinor bundle on $Q_{2m+1}$. Thus we have that
if $n\geq 5$ then $H^1_*(F\otimes\Sigma )=0$ because $F$ is ACM and $H^2_*(F\otimes\Sigma)=0$.

If $n=3$, using the restriction sequence for the bundle $E\otimes\Sigma$ and (iv), we get that
$H^1(E\otimes\Sigma(k))=0$ for $k\le-1$ and $H^2(E\otimes\Sigma(k))=0$ for $k\ge-1$.
Then using the two short exact sequences as in (\ref{equker}) and (\ref{equimg}) tensored
by the spinor bundle $\Sigma$, we get that $H^1(F\otimes\Sigma(k))=0$ for $k\le-1$ and
$H^2(F\otimes\Sigma(k))=0$ for $k\ge-1$. Finally, an exact sequence like (\ref{fsigma})
allows us to conclude that $H^1(F\otimes\Sigma(k))=0$ for $k\ge-1$ as well, therefore
$H^1_*(F\otimes\Sigma)=0$.

Similarly, for $n=2m$ ($m\ge2$) tensoring the exact sequences of
\cite[Theorem 2.8(ii)]{O1} by $F(k)$ we have:
$$
0 \to F\otimes\Sigma_1 (k)\to F(k)^{\oplus 2^{m}} \to F\otimes\Sigma_2(k+1) \to 0
$$
and
$$
0 \to F\otimes\Sigma_2 (k) \to F(k)^{\oplus 2^{m}} \to F\otimes\Sigma_1(k+1) \to 0
$$
where $\Sigma_1$ and $\Sigma_2$ are the two spinor bundles on $Q_{2m}$. Thus we also have that
$H^1_*(F\otimes\Sigma _{i})=0$ for $i=1,2$.

>From \cite[Theorem 3.4]{Ot} we have that a bundle $F$ over $Q_{n}$ splits
as a sum of line bundles if and only if $H^p_*(F)=0$ for $2\le p\le n-1$ and
$H^p_*(F\otimes\Sigma)=0$ for every spinor bundle $\Sigma$ and each $1\le p\le n-2$. Hence
we can write $F=\oplus _{i=1}^{r}\oo _{Q_{n}}(k_{i})$. 

Moreover, from (\ref{equofo}) we have that  ${\rm ch}(E)={\rm ch}(F)-a'-a\cdot{\rm ch}(\oo _{Q_{n}}(-1))-c'-
c\cdot {\rm ch}(\oo _{Q_{n}}(1))$ and using (i) we are led to ${\rm ch}(F)=a'+b+c'$.
Since we have that 
${\rm ch}(F)=r+(\sum _{i=1}^{r}k_{i})h+(\sum _{i=1}^{r}k_{i}^{2})h^{2}+\ldots$,
it follows that $r=a'+b+c'$ and $\sum _{i=1}^{r}k_{i}^{2}=0$. Hence $E$ is the
cohomology of a monad of the form
$$ 0\longrightarrow\oo _{Q_n}(-1)^{\oplus a}\oplus\oo _{Q_n}^{\oplus a'}\stackrel{\alpha}\longrightarrow
\oo _{Q_n}^{\oplus a'+b+c'}\stackrel{\beta}
\longrightarrow\oo _{Q_n}^{\oplus c'}\oplus\oo _{Q_{n}}(1)^{\oplus c}\longrightarrow 0. $$

Now consider the restricted morphism $\alpha :\oo _{Q_{n}}\rightarrow\oo _{Q_{n}}^{\oplus r}$ where 
$r=a'+b+c'$. We have that $\alpha$ corresponds to a morphism of graded algebras 
$\alpha ':S(Q_{n})\rightarrow S(Q_{n})^{\oplus r}$ which is completely determined by the image 
of $1\in S(Q_{n})_{0}$. Therefore, reordering the factors of $S(Q_{n})^{\oplus r}$ if necessary,
we can write $\alpha '(1)=(1,\ldots ,1,0,\ldots ,0)$ and so 
$\alpha '(S(Q_{n}))=(1,\ldots ,1,0,\ldots ,0)\ S(Q_{n})$ what implies $S(Q_{n})^{\oplus r}/\alpha '(S(Q_{n}))\cong
S(Q_{n})^{r-1}$ which yields $\oo _{Q_{n}}^{\oplus r}/\alpha (\oo _{Q_{n}})\cong\oo _{Q_{n}}^{\oplus r-1}$.

Since $E$ remains the cohomology of the monad
$$ 0\longrightarrow\oo _{Q_n}(-1)^{\oplus a}\oplus\oo _{Q_n}^{\oplus a'}/\oo _{Q_{n}}\stackrel{\bar{\alpha}}\longrightarrow
\oo _{Q_n}^{\oplus a'+b+c'}/\alpha (\oo _{Q_{n}})\stackrel{\beta}
\longrightarrow\oo _{Q_n}^{\oplus c'}\oplus\oo _{Q_{n}}(1)^{\oplus c}\longrightarrow 0$$
we can repeat this procedure to cancel the summands $\oo _{Q_{n}}^{\oplus a'}$'s from first and middle terms.
In order to cancel the $\oo _{Q_{n}}^{\oplus c'}$'s from middle and last terms we can pass to the dual monad, 
reapply the same argument above to cancel them and get back to the original monad dualizing once more. 
The resulting monad will fulfill the desired form. 
\end{proof}

It follows from Theorem \ref{thmhyp} that any extension of linear bundles 
on $Q_n$ ($n\ge3$) is also a linear bundle. As in the case of projective spaces, linear
bundles on $Q_n$ exist if and only if at least either $b\geq 2c+n-1$ and 
$b\geq a+c$ or $b\geq a+c+n$, see \cite[Proposition 4.7]{CMR}.


\section{Steiner Bundles}
\label{secstn}

Let $X$ be a projective variety. We will say that a sheaf $S$ on $X$ is {\em Steiner}
if it has a resolution of the form:
\begin{equation}
\label{steiner}
0\longrightarrow\oo _{X}(-1)^{\oplus a}\longrightarrow
\oo _{X}^{\oplus b}\longrightarrow S \longrightarrow 0.
\end{equation}
As it was mentioned at the Introduction, we are just extending in a natural way the
definition introduced by I. Dolgachev and M. Kapranov for the case of projective spaces
in \cite{DK}, where Steiner bundles appeared in the study of hyperplane arrangements.
There they used the Beilinson spectral sequence to establish the following cohomological characterization (cf. \cite[Proposition 3.2]{DK}): a locally free sheaf $S\to\pp^n$ is Steiner if and only if $H^q(S\otimes\Omega^p_{\pp ^n}(p)) = 0$ for $q>0$ and for $q=0,~p>1$. Steiner bundles exist provided ${\rm rk}(S)=b-a\ge n$; in the critical case, when $n=b-a$, Steiner bundles are also known as Schwarzenberger bundles (cf. \cite{AO}).

\subsection{Steiner bundles on Projective Spaces}

Basically with the same technique used before we will prove here the following alternative cohomological characterization
of Steiner bundles:

\begin{teo}
\label{thmstp}
A locally free sheaf $S$ on $\pp^n$ with $n\geq 3$ is Steiner if and only if:
\begin{itemize}
\item[\um] ${\rm ch}(S)=b-a\cdot{\rm ch}(\oo _{\pp ^{n}}(-1))$;
\item[\doi] $H^0(S(-1))=H^{n-1}(S(1-n))=H^n(S(-n))=0$;
\item[\tre] $H^p_*(S)=0$ where $1\le p\le n-2$.
\end{itemize}
\end{teo}

\begin{proof}
The ``only if" part is a straightforward calculation starting with the exact sequence
$$ 0\longrightarrow\oo _{\pp ^n}(k-1)^{\oplus a}\longrightarrow
\oo _{\pp ^{n}}^{\oplus b}(k)\longrightarrow S(k) \longrightarrow 0.$$
Conversely, note that (ii) and (iii) imply via Serre duality:
\begin{itemize}
\item[(ii')] $H^0(S^*(-1))=H^{1}(S^*(-2))=H^n(S^*(-n))=0$;
\item[(iii')] $H^p_*(S^{*})=0$ where $2\le p\le n-1$.
\end{itemize}
It then follows, in particular, that $H^{1}(S^*(k))=0$ for $k\le-2$; hence by Lemma \ref{lemgen}
with $l=-1$, we conclude that $H^{1}(S^*(-1))$ generates $H^{1}_*(S^*)$ as $S(\pp ^n)$-module.

Now let us take a basis $g_{1},\ldots ,g_{a}$ of $H^1(S^*(-1))$ where $a=h^1(S^*(-1))$.
The element $\sum _{i} g_{i}\in H^{1}(S^{*}(-1))^{\oplus a}\cong
\text{Ext}^{1}(\oo _{\pp ^{n}}(1)^{\oplus a},S^{*})$ yields an extension
$$ 0\longrightarrow S^* \longrightarrow F \longrightarrow
\oo _{\pp ^n}(-1)^{\oplus a} \longrightarrow 0.$$
>From (iii') we have that $H^{p}_{*}(F)=0$ for $2\le p\le n-1$. And from the proof of Theorem \ref{thmhor}, 
we have that $H^{1}_{*}(F)=0$ if and only if the $g_{i}'s$ generate $H^{1}_{*}(S^{*})$ as $S(\pp ^n)$-module. 
Since we have just seen that this is the case it follows that $F$ is ACM. 

Thus by Horrocks' splitting criterion
$F=\oplus_{i=1}^{b} \oo_{\pp^n}(a_i)$. Moreover, (i) implies $ch(F)=b$ and so
${\rm ch}_2(F)=\frac{1}{2}\sum_{i=1}^{b} a_i^2=0$. Hence the 
$a_i$'s must vanish and therefore $F=\oo_{\pp^n}^{\oplus b}$. Dualizing the above exact sequence we are led
to the desired resolution.
\end{proof}

It follows from Theorem \ref{thmstp} that any extension of Steiner bundles on 
$\pp^n$ is also a Steiner bundle.

\subsection{Steiner bundles on quadric hypersurfaces}

Now we give a cohomological characterization of Steiner bundles on quadric hypersurfaces. 

\begin{teo}
\label{thmstq}
A locally free sheaf $S$ on $Q_n$ with $n\ge3$ is Steiner if and only if:
\begin{itemize}
\item[\um] ${\rm ch}(S)=b-a\cdot{\rm ch}(\oo _{Q_n}(-1))$;
\item[\doi] $H^0(S(-1))=H^{n-1}(S(2-n))=H^n(S(1-n))=0$;
\item[\tre] $H^p_*(S)=0$ where $1\le p\le n-2$;
\item[\qua] $H^p_*(S\otimes\Sigma)=0$ for every spinor bundle $\Sigma$ and each $1\le p\le n-2$.
\end{itemize}
\end{teo}

\begin{proof}
The ``only if" part is a straightforward calculation starting with the exact sequence
$$ 0\longrightarrow\oo _{Q_n}(k-1)^{\oplus a}\longrightarrow
\oo _{Q_n}^{\oplus b}(k)\longrightarrow S(k) \longrightarrow 0.$$
Conversely, hypotheses (ii) and (iii) 
imply via Serre duality:
\begin{itemize}
\item[(ii')] for $n\ge2$, $H^0(S^*(-1))=H^{1}(S^*(-2))=H^n(S^*(1-n))=0$;
\item[(iii')] for $n\ge3$, $H^p_*(S^*)=0$ where $2\le p\le n-1$.
\end{itemize}
It then follows that $H^{1}(S^*(k))=0$ for $k\le-2$; hence by Lemma \ref{lemgen}
with $l=0$, we conclude that $H^{1}(S^*(-1))\oplus H^{1}(S^*)$
generates $H^{1}_{*}(S^*)$ as $S(Q_n)$-module. Then following the proof of Theorem
\ref{thmstp}, we have the exact sequence:
\begin{equation}\label{dualsqc}
0\longrightarrow S^* \longrightarrow F \longrightarrow
\oo _{Q_n}(1)^{\oplus a}\oplus\oo _{Q_n}^{\oplus a'}\longrightarrow 0,
\end{equation}
with $a=h^1(S^*(-1))$ and $F$ being ACM. Now (iv) implies that
$H^p_*(F\otimes\Sigma)=0$ for $1\le p\le n-2$, hence $F$ must split as a sum
of line bundles by \cite[Theorem 3.4]{Ot}.
Then (i) forces $ch(F)=b+a'$ and hence $F=\oo_{Q_n}^{\oplus b+a'}$. Dualizing
sequence (\ref{dualsqc}) we obtain:
$$0\longrightarrow \oo _{Q_n}(-1)^{\oplus a}\oplus\oo _{Q_n}^{\oplus a'}
\longrightarrow \oo _{Q_n}^{\oplus b+a'} \longrightarrow S \longrightarrow 0$$
and we apply the same argument at the end of the proof of Theorem \ref{thmhyp} to cancel the
extra summand $\oo _{Q_n}^{\oplus a'}$.
\end{proof}

As in the case of projective spaces, it follows from Theorem \ref{thmstq} that
any extension of Steiner bundles on $Q_n$ is also a Steiner bundle.


\section{Linear resolutions}
\label{secres}

Let $E$ be a sheaf on a projective variety $X$. Following \cite{BS},
a linear resolution of $E$ is an exact sequence
$$ 0 \to F_d \to F_{d-1} \to \cdots \to F_1 \to F_0 \to E \to 0 $$
where each $F_k$ splits as a sum of line bundles; we call the minimal number of terms
in such resolutions the {\em linear homological dimension} of $E$, which is denoted by
${\rm lhd}(E)$.

\begin{prop}
\label{prpres}
Let $E$ be a locally free sheaf on an ACM variety $X$. We have:
$$
{\rm lhd}(E)\le d\ \ \Longrightarrow\ \ H^p_*(E\otimes W)=0
$$
for every ACM locally free sheaf $W$ on $X$ and each $1\le p\le n-d-1$.
\end{prop}

\begin{proof}
We proceed by induction on $d$. For $d=0$, the statement is simply
the definition of an ACM locally-free sheaf. So assume that if ${\rm lhd}(K)\le d-1$,
then $H^p_*(K\otimes W)=0$ for every ACM locally-free sheaf $W$ and each $1\le p\le n-d$.
If $E$ has ${\rm lhd}(E)\le d$, then there is an exact sequence
$$ 0 \to F_d \to F_{d-1} \to \cdots \to F_1 \to F_0
\stackrel{\phi}{\longrightarrow} E \to 0 $$
where each $F_k$ splits as a sum of line bundles. Taking $K=\ker\phi$, we have two
exact sequences:
$$ 0 \to F_d \to F_{d-1} \to \cdots \to F_1 \to K \to 0, $$
which means that ${\rm lhd}(K)\le d-1$, and
\begin{equation}
\label{kerphi}
0 \to K \to F_0 \stackrel{\phi}{\longrightarrow} E \to 0.
\end{equation}
The Proposition is then easily proved after tensoring the last sequence with $W$ and
applying the induction hypothesis.
\end{proof}

For $X=\pp^n$, it was proved by Bohnhorst and Spindler that the converse of
Proposition \ref{prpres} is also true, see \cite[Proposition 1.4]{BS}:
if $E$ is a locally-free sheaf on $\pp^n$ satisfying $H^p_*(E)$ for $1\le p\le n-d-1$,
then ${\rm lhd}(E)\le d$. Let us now establish the converse of Proposition \ref{prpres}
to smooth quadric hypersurfaces, following \cite{BS}. We will require the following
technical result.

\begin{lema}
\label{surjh0}
Let $E$ be a locally free sheaf on a projective variety $X$. Then there are
$a_1,\dots,a_k$ and a surjective map $\phi:\oplus_{i=1}^k\oox(-a_i)\to E$
such that for each integer $m$ the induced map in cohomology:
$$ H^0(\oplus_{i=1}^k\oox(m-a_i)) \to H^0(E(m)) $$
is also surjective.
\end{lema}
\begin{proof}
Since $E$ is locally free, there is an integer $m_0$ such that
$H^0(E(m))=0$ for each $m\le m_0$, and $H^0_*(E)$ is a finitely
generated graded $S(X)$-module. Thus there are integers
$a_1,\dots,a_k$ (the degrees of generators $g_1\dots,g_k$ of
$H^0_*(E)$) for which the map
\begin{gather*}
\begin{matrix}
\oplus_{i=1}^k H^0(\oox(m-a_i)) & \longrightarrow & H^0(E(m))\\
(f_1,\dots,f_k) & \longmapsto & \sum_{i=1}^n f_ig_i
\end{matrix}
\end{gather*}
is surjective for each $m$. But this is exactly the induced map
in cohomology of a surjective sheaf map $\oplus_{i=1}^k\oox(-a_i)\to E$,
as desired.
\end{proof}

Notice that the sheaf $K=\ker\phi$ is such that $H^1_*(K)=0$ and
$H^{p+1}_*(K)\cong H^{p}_*(E)$ for $1\le p\le n-2$. Next, we prove
a generalization of \cite[Proposition 1.4]{BS} to general ACM varieties

\begin{prop}
\label{longacm}
Let $E$ be a locally free sheaf on an ACM projective variety $X$.
$H^p_*(E)=0$ for $1\le p\le n-d-1$ if and only if there is an exact sequence
$$ 0 \to K_d \to F_{d-1} \to \cdots \to F_1 \to F_0 \to E \to 0 $$
where each $F_k$ splits as a sum of line bundles and $K_d$ is an ACM
locally free sheaf. 
\end{prop}
\begin{proof}
Again, the proof is by induction on $d\ge1$. From the previous lemma 
we know that, for any locally free sheaf $E$, there is an exact sequence
\begin{equation}
\label{equkph}
0\to K_1 \to L_0 \stackrel{\phi}{\to} E \to 0
\end{equation}
such that $L_0$ splits as a sum of line bundles and $H^1_*(K_1)=0$. 
If moreover $H^p_*(E)=0$ for $1\le p\le n-2$, then also $H^p_*(K_1)=0$ for
$2\le p\le n-1$, hence $K_1$ is ACM; this concludes the proof of the case
$d=1$.

Now assume that $H^p_*(E)=0$ for $1\le p\le n-d-1$. Following the argument of the
previous paragraph, we find an exact sequence
$$ 0\to K_1 \to L_0 \stackrel{\phi}{\to} E \to 0 $$
such that $H^p_*(K_1)=0$ for $1\le p\le n-d$. Thus by the induction hypothesis,
there exists a long exact sequence
$$ 0 \to K_{d-1} \to F_{d-2} \to \cdots \to F_1 \to F_0 \to K_1 \to 0 $$
where each $F_k$ splits as a sum of line bundles and $K_{d-1}$ is an ACM
locally-free sheaf. Glueing the two last exacts sequences, we obtain
$$ 0 \to K_{d-1} \to F_{d-2} \to \cdots \to F_0 \to L_0 \to E \to 0$$
as desired.

The converse statement is also deduced by induction on $d$.
If there is an exact sequence
$$ 0\to K_1 \to F_0 \stackrel{\phi}{\to} E \to 0 $$
with $K_1$ being ACM and $F_0$ splitting as a sum of line bundles,
it follows that $H^{p}_*(E)\cong H^{p+1}_*(K_1)=0$ for $1\le p\le n-2$,
which proves the case $d=1$. Now assume there is an exact sequence
$$ 0 \to K_d \to F_{d-1} \to \cdots \to F_1 \to F_0 \to E \to 0 $$
where each $F_k$ splits as a sum of line bundles and $K_d$ is ACM.
It can be broken into two exact sequences:
$$ 0 \to K_d \to F_{d-1} \to \cdots \to F_1 \to K_1 \to 0 $$
and 
$$ 0\to K_1 \to F_0 \stackrel{\phi}{\to} E \to 0 ~~. $$
By the first sequence and the induction hypothesis, we conclude that
$H^p_*(K_1)=0$ for $1\leq p \leq n-(d-1)-1=n-d$. Then use the
second sequence to deduce that $H^p_*(E)=0$ for $1\leq p \leq n-d-1$.
\end{proof}

Now we are finally in a position to prove our last result.

\begin{teo}
\label{thmres}
Let $E$ be a locally free sheaf on $Q_n$ with $n\ge 3$. We have:
$$ {\rm lhd}(E)\le d ~~\Longleftrightarrow~~
\left\{ \begin{array}{l} H^p_*(E)=0\\
H^{n-d-1}_*(E\otimes\Sigma)=0
\end{array} \right.$$
for every spinor bundle $\Sigma$ and each $1\leq p\leq n-d-1$.
\end{teo}

\begin{proof}
It is enough to prove the ``only if" direction. From Proposition \ref{longacm} we know
that since $H^p_*(E)=0$ for $1\le p\le n-d-1$ there is a long exact sequence  
$$ 0 \to K_d \to F_{d-1} \to \cdots \to F_1 \to F_0 \to E \to 0 $$
where each $F_k$ splits as a sum of line bundles and $K_d$ is an ACM bundle. Breaking it
into short exact sequences and passing to cohomology, one can see that 
$H^{n-d-1}_*(E\otimes\Sigma)=0$ implies that $H^{n-1}_*(K_d\otimes\Sigma)=0$.
But now from \cite[Corollary 4.3]{CMR3} we have that a locally free sheaf $E$ 
on $Q_n$ ($n\ge 3$) splits as a sum of line bundles if and only if it is ACM and $H^{n-1}_*(E\otimes\Sigma)=0$ 
for every spinor bundle $\Sigma$ on $Q_n$. It then follows that $K_d$ must split as a sum of line bundles and
hence ${\rm lhd}(E)\le d$ as desired.
\end{proof}

\end{document}